\documentstyle[12pt]{article}

\begin{document}
\newcommand{\ol }{\overline}
\newcommand{\ul }{\underline }
\newcommand{\ra }{\rightarrow }
\newcommand{\lra }{\longrightarrow }
\newcommand{\ga }{\gamma }
\newcommand{\st }{\stackrel }
\newcommand{\scr }{\scriptsize }
\title{{\Large\bf On Polynilpotent Covering Groups
of a Polynilpotent Group}
\vspace{-0.7cm}\author{ Behrooz Mashayekhy and Mahboobeh Alizadeh Sanati \\
Department of Mathematics, Ferdowsi University of Mashhad,\\
P.O.Box 1159-91775, Mashhad, Iran\\
E-mail:  mashaf@math.um.ac.ir\\
\hspace{2.05cm}sanati@math.um.ac.ir} }
\date{ }
\maketitle
\vspace{-0.7cm}
\begin{abstract}
 Let ${\cal N}_{c_1,\ldots,c_t}$ be the variety of polynilpotent
 groups of class row $(c_1,\ldots,c_t)$. In this paper, first, we
 show that a polynilpotent group $G$ of class row
 $(c_1,\ldots,c_t)$ has no any ${\cal N}_{c_1,\ldots,c_t,c_{t+1}}$-covering
 group if its Baer-invariant with respect to the variety  ${\cal
 N}_{c_1,\ldots,c_t,c_{t+1}}$ is nontrivial. As an immediate
 consequence, we can conclude that a solvable group $G$ of length
 $c$ with nontrivial solvable multiplier, ${\cal S}_nM(G)$, has no
 ${\cal S}_n$-covering group for all $n>c$, where ${\cal S}_n$ is
 the variety of solvable groups of length at most $n$. Second, we
 prove that if $G$ is a polynilpotent group of class row
 $(c_1,\ldots,c_t,c_{t+1})$ such that ${\cal N}_{c'_1,\ldots,c'_t,c'_{t+1}}M(G)\neq
 1$, where $c'_i\geq c_i$ for all $1\leq i\leq t$ and
 $c'_{t+1}>c_{t+1}$, then $G$ has no any ${\cal
 N}_{c'_1,\ldots,c'_t,c'_{t+1}}$-covering group. This is a vast
 generalization of the first author's result on nilpotent covering
 groups (Indian\ J.\ Pure\ Appl.\ Math.\ 29(7)\ 711-713,\ 1998).



\end{abstract}
\begin{center}
{\bf 1.Introduction and Motivation}
\end{center}
{\hspace{0.6cm}}Let $G\simeq F/R$ be a free presentation for $G$ and $\cal V$ be a
variety of groups. Then, after R. Baer [1], the Baer-invariant of
$G$ with respect to $\cal V$ is defined to be ${\cal V}M(G)=R\cap
V(F)/[RV^*F]$, where $V(F)$ is the verbal subgroup of $F$ with
respect to $\cal V$ and\\
$[RV^*F]=\langle
v(f_1,\cdots,f_{i-1},f_ir,f_{i+1},\cdots,f_n)v(f_1,\cdots,f_n)^{-1}\
|~ r\in R,~ f_i\in F,$
$$1\leq i\leq n,~ v\in V, ~n\in{\bf N}\rangle\ .$$
{\hspace{0.6cm}}In special case, if $\cal V$ is the variety of abelian groups,
then the Baer-invariant of $G$ will be the well-known notion the
Schur-multiplier of $G$, denoted by $M(G)=R\cap F'/[R,F]$ ( See
[5,6] for further details).

It is easy to see that if ${\cal V=N}_c$, the variety of
nilpotent groups of class at most $c\geq 1$, then
$${\cal N}_cM(G)=\frac {R\cap \gamma_{c+1}(F)}{[R, _cF]}\ \ ,$$
where $\gamma_{c+1}(F)$ is the $(c+1)$-st term of the lower central
series of $F$ and $[R,_1F]=[R,F]$, $[R,_cF]=[[R,_{c-1}F],F]$,
inductively. We shall also call ${\cal N}_cM(G)$ the
$c$-nilpotent multiplier of $G$.

In a more general case, if ${\cal V}={\cal
 N}_{c_1,\ldots,c_t}$, the variety of polynilpotent groups of
 class row $(c_1,\ldots,c_t)$, then
$${\cal N}_{c_1,\ldots,c_t}M(G)=\frac{R\cap\gamma_{c_t+1}\circ\cdots\circ\gamma_{c_1+1}(F)}{[R,_{c_1}F,_{c_2}\gamma_{c_1+1}(F),\cdots,_{c_t}\gamma_{c_{t-1}+1}\circ\cdots\circ\gamma_{c_1+1}(F)]}\ \ ,$$
where
$\gamma_{c_t+1}\circ\cdots\circ\gamma_{c_1+1}(F)=\gamma_{c_t+1}(\gamma_{c_{t-1}+1}(\cdots(\gamma_{c_1+1}(F))\cdots))$
are the terms of iterated lower central series of $F$. See [4,
corollary 6.14] for the following equality
$$[RN^*_{c_1,\ldots,c_t}F]=[R,_{c_1}F,_{c_2}\gamma_{c_1+1}(F),\ldots,_{c_t}\gamma_{c_{t-1}+1}\circ\ldots\circ\gamma_{c_1+1}(F)].$$
We shall also call ${\cal N}_{c_1,\ldots,c_t}M(G)$, the
$(c_1,\ldots,c_t)$-polynilpotent multiplier of $G$.

Let $\cal V$ be a variety of groups and $G$ be an arbitrary
group, then a $\cal V$-covering group of $G$ (a generalized
covering group of $G$ with respect to the variety $\cal V$) is a group $G^*$
with a normal subgroup $A$ such that $G^*/A\simeq G$, $A\subseteq
V(G^*)\cap V^*(G^*)$, and $A\simeq{\cal V}M(G)$, where $V^*(G^*)$
is the marginal subgroup of $G^*$ with respect to $\cal V$ (see
[6]).

Note that if $\cal V$ is the variety of abelian groups, then the
$\cal V$-covering group of $G$ will be ordinary covering group (sometimes it is called representing group) of $G$. Also if ${\cal
V}={\cal N}_{c_1,\ldots,c_t}$, then an ${\cal
 N}_{c_1,\ldots,c_t}$-covering group of $G$ is a group $G^*$ with
 a normal subgroup $A$ such that

 $G\simeq G^*/A$,

 $A\simeq{\cal N}_{c_1,\ldots,c_t}M(G^*)$ and

 $A\subseteq
 N^*_{c_1,\ldots,c_t}(G^*)\cap\gamma_{c_t+1}(\cdots(\gamma_{c_1+1}(G^*))\cdots)$.\\
We shall also call $G^*$ a $(c_1,\ldots,c_t)$-polynilpotent
covering group of $G$.

It is a well-known fact that every group has at least a covering
group (see [5,13]). Also, the first author proved that every group
has a ${\cal V}$-covering group if ${\cal V}$ is the variety of
all groups, $\cal G$, or the variety of all abelian groups, ${\cal
A}$, or the variety of all abelian groups of exponent $m, {\cal
A}_{m}$, where $m$ is square free (see [7,9]).

Moreover, C. R. Leedham-Green and S. Mckay [6] proved, by a
homological method, that a sufficient condition for existence of
a ${\cal V}$-covering group of $G$ is that $G/V(G)$ should be a
${\cal V}$-splitting group.

Some people have tried to construct a covering group for some
well-known structures of groups. For example, the generalized
quaternion  group $Q_{4n}=\langle a, b | a^{2n}=1, b^{2}=a^{n},
b^{-1}ab=a^{-1}\rangle$ is a covering group of the dihedral group
$D_{2n}=\langle a, b | a^n=b^2=1, b^{-1}ab=a^{-1}\rangle$ (see
[5]).

Also J. Wiegold [12] presented a covering group for a direct
product of two finite groups. W. Haebich [2, 3] generalized the
Wiegold's result and gave a covering group for a regular product
of a family of groups and also for a verbal wreath product of two
groups. Moreover, the first author [10] recently proved the
existence and presented a structure of an ${\cal N}_c$-covering
group for a nilpotent product of a family of cyclic groups.

It is interesting to mention that there are some groups which have
no any $\cal V$-covering group, for some variety $\cal V$. The
first author [8] gave an example of the group $G\simeq{\bf
Z}_r\oplus{\bf Z}_s$, where $(r,s)\neq 1$, which has no ${\cal
N}_c$-covering group for all $c\geq 2$. Moreover, the first author
[7, 9] proved that a nilpotent group $G$ of class $n$ with
nontrivial $c$-nilpotent multiplier ${\cal N}_cM(G)$, has no
${\cal N}_c$-covering group for all $c>n$.

Now, in this paper, we concentrate on nonexistence of polynilpotent
covering groups of a polynilpotent group. More precisely, we show
that if $G$ is a polynilpotent group of class row
$(c_1,\ldots,c_t)$ such that ${\cal
N}_{c_1,\ldots,c_t,c_{t+1}}M(G)\neq 1$, then $G$ has no
$(c_1,\ldots,c_t,c_{t+1})$-polynilpotent covering group of $G$.
Also, if ${\cal N}_{c'_1,\ldots,c'_t}M(G)\neq 1$ and $c'_i\geq
c_i$ for all $1\leq i\leq t-1$ and $c'_t>c_t$, then $G$ has no
$(c'_1,\ldots,c'_t)$-polynilpotent covering group of $G$.
\begin{center}
{\bf 2.The Main Results}
\end{center}
{\hspace{0.6cm}}Let $G$ be a group and $\cal V$ be a variety of groups. It is
clear, by definition, that if ${\cal V}M(G)=1$, then $G$ is the
only $\cal V$-covering group of itself. So it is natural to put
the condition ${\cal V}M(G)\neq 1$ for nonexistence of $\cal
V$-covering group of $G$.\\
{\bf Theorem 2.1}

Let $G$ be a polynilpotent group of class row $(c_1,\ldots,c_t)$
and ${\cal N}_{c_1,\ldots,c_t,c_{t+1}}M(G)$\\
$\neq 1$, for some $c_{t+1}\geq 1$. Then $G$ has no any ${\cal N
}_{c_1,\ldots,c_t,c_{t+1}}$-covering group.\\
${\bf Proof.}$

Let $G^*$ be a $(c_1,\ldots,c_t,c_{t+1})$-polynilpotent covering
group of $G$ with the normal subgroup $A$ of $G^*$ such that

$G\simeq G^*/A$,

$A\simeq{\cal N}_{c_1,\ldots,c_t,c_{t+1}}M(G^*)$ and

$A\subseteq N^*_{c_1,\ldots,c_t,c_{t+1}}(G^*)\cap\gamma_{c_{t+1}+1}(\gamma_{c_t+1}(\cdots(\gamma_{c_1+1}(G^*))\cdots))$.\\
We define $\rho_t$ inductively, for any group $M$ and $t\geq 0$,
as follows:
\begin{center}
$\rho_0(M)=M$ and $\rho_i(M)=\gamma_{c_i+1}(\rho_{i-1}(M)),$ for
$i>1$.
\end{center}
By hypothesis, $\rho_t(G)=1$ and so $\rho_t(G^*/A)=1$. Hence
$\rho_t(G^*)\subseteq A$. Also $A\subseteq\rho_{t+1}(G^*),$ then
$\rho_t(G^*)\subseteq\rho_{t+1}(G^*)$. Clearly
$\rho_{t+1}(G^*)\subseteq\rho_t(G^*)$, so
$\rho_{t+1}(G^*)=\rho_t(G^*)$. In particular,\\
$\rho_t(G^*)=\gamma_2(\rho_t(G^*))=\cdots=\gamma_{c_{t+1}}(\rho_t(G^*))=\gamma_{c_{t+1}+1}(\rho_t(G^*))=\rho_{t+1}(G^*)\
(I).$\\
 Since $A\subseteq N^*_{c_1,\ldots,c_t,c_{t+1}}(G^*)$ and
$\rho_t(G^*)\subseteq A$, so we have
$$ [\cdots[[\rho_t(G^*),_{c_1}G^*],_{c_2}\gamma_{c_1+1}(G^*)],\cdots,_{c_{t+1}}\gamma_{c_t+1}(\cdots(\gamma_{c_1+1}(G^*))\cdots)]=1,$$
or by the above notation,
$$[\cdots[[\rho_t(G^*),_{c_1}G^*],_{c_2}\rho_1(G^*)],\cdots,_{c_{t+1}}\rho_t(G^*)]=1.$$
{\hspace{0.6cm}}First, we show that $[M,_iN]\stackrel{(II)}{\supseteq}[\gamma_i(N),M]$ for each
natural number i and normal subgroups $M$ and $N$ of any group. By Three Subgroups Lemma, we have\\
$[M,_iN]=[M,_{i-2}N,N,N]\supseteq[N,N,[M,_{i-2}N]]=[[M,_{i-2}N,[N,N]]$\\
\hspace{1cm}$=[[M,_{i-3}N],N,\gamma_2(N)]\supseteq[N,\gamma_2(N),[M,_{i-3}N]]=[[M,_{i-3}N],\gamma_3(N)]$\\
\hspace{.5cm}$=[[M,_{i-4}N],N,\gamma_3(N)]\supseteq\cdots\supseteq[M,\gamma_i(N)]=[\gamma_i(N),M].$

Now, we claim\\
$\hspace{1.4cm}[\cdots[[\rho_t(G^*),_{c_1}G^*],_{c_2}\rho_1(G^*)],\cdots,_{c_i}\rho_{i-1}(G^*)]\stackrel{(III)}{\supseteq}[\rho_i(G^*),\rho_t(G^*)]\
,$\\
for all $1\leq i\leq t+1$.

Clearly the equality is valid for $i=1$. Now for $i=2$, we can write\\
$[[\rho_t(G^*),_{c_1}G^*],_{c_2}\rho_1(G^*)]\supseteq[[\gamma_{c_1}(G^*),\rho_t(G^*)],_{c_2}\rho_1(G^*)]\hspace{0.8cm}by (II)$\\
$\supseteq[[\rho_1(G^*),\rho_t(G^*)],_{c_2}\rho_1(G^*)]$\\
$\supseteq[\gamma_{c_2}(\rho_1(G^*)),[\rho_1(G^*),\rho_t(G^*)]]\hspace{1.8cm}by (II)$\\
$=[[\rho_1(G^*),\rho_t(G^*)],\gamma_{c_2}(\rho_1(G^*))]$\\
$\supseteq[\gamma_{c_2}(\rho_1(G^*)),\rho_1(G^*),\rho_t(G^*)]$\\
$=[\gamma_{c_2+1}(\rho_1(G^*)),\rho_t(G^*)]$\\
$=[\rho_2(G^*),\rho_t(G^*)]$.\\
Suppose the inclusion (III) holds for $i=j$. Now, we prove it for $i=j+1$.\\
$[[\cdots[[\rho_t(G^*),_{c_1}\rho_0(G^*)],_{c_2}\rho_1(G^*)],\cdots,_{c_j}\rho_{j-1}(G^*)],_{c_{j+1}}\rho_j(G^*)]$\\
$\supseteq[\rho_j(G^*),\rho_t(G^*),_{c_{j+1}}\rho_j(G^*)]$\\
$\supseteq[\gamma_{c_{j+1}}(\rho_j(G^*)),[\rho_j(G^*),\rho_t(G^*)]]\hspace{3cm}by (I)$\\
$=[\rho_j(G^*),\rho_t(G^*),\gamma_{c_{j+1}}(\rho_j(G^*))]$\\
$\supseteq[\gamma_{c_{j+1}}(\rho_j(G^*)),\rho_j(G^*),\rho_t(G^*)]$\\
$=[\gamma_{c_{j+1}+1}(\rho_j(G^*)),\rho_t(G^*)]$\\
$=[\rho_{j+1}(G^*),\rho_t(G^*)].$\\
Now, we have
$$1=[\cdots[[\rho_t(G^*),_{c_1}\rho_0(G^*)],_{c_2}\rho_1(G^*)],\cdots,_{c_{t+1}}\rho_t(G^*)]\supseteq[\rho_{t+1}(G^*),\rho_t(G^*)].$$
Hence $[\rho_{t+1}(G^*),\rho_t(G^*)]=1$. Since
$\rho_{t+1}(G^*)=\rho_t(G^*)$, we can conclude
$[\rho_{t}(G^*),\rho_t(G^*)]=1$. i.e. $\gamma_2(\rho_t(G^*))=1$.
Hence by $(I)$, we have $\rho_{t+1}(G^*)=1$. Therefore $A=1$,
which is a contradiction. $\Box$

Now we can state the following interesting corollary about
nonexistence of solvable covering groups.\\
{\bf Corollary 2.2}

Let $G$ be a solvable group with derived length at most $n$. If
the $l$-solvable multiplier of $G$, ${\cal S}_lM(G)$, is
nontrivial, then $G$ has no any ${\cal S}_l$-covering group, for
all $l>n$.\\
${\bf Proof.}$

Note that, ${\cal S}_l$, the variety of solvable groups of
derived length at most $l$ is in fact the variety of
polynilpotent groups of class row
$\underbrace{(1,\ldots,1)}_{l-times}$. Hence the result is a
consequence of Theorem 2.1. $\Box$\\

In a different view, the following theorem is also about
nonexistence of polynilpotent covering groups which is a vast
generalization of a result of the first author (see [7, Theorem
3.1.6], [8,Theorem 2] and [9, Theorem 2.1]).\\
{\bf Theorem 2.3}

Let $G$ be a polynilpotent group of class row
$(c_1,\ldots,c_t,c_{t+1})$ such that ${\cal
N}_{c'_1,\ldots,c'_t,c'_{t+1}}M(G)\neq 1$ where $c'_i\geq c_i$ for all $1\leq i\leq t$ and $c'_{t+1}>c_{t+1}$. Then $G$ has no any
${\cal
N}_{c'_1,\ldots,c'_t,c'_{t+1}}$-covering group.\\
${\bf Proof.}$

Let $G^*$ be a $(c'_1,\ldots,c'_t,c'_{t+1})$-polynilpotent covering
group of $G$ with the normal subgroup $A$ of $G^*$ such that

$G\simeq G^*/A$,

$A\simeq{\cal N}_{c'_1,\ldots,c'_t,c'_{t+1}}M^*(G)$ and

$A\subseteq
N^*_{c'_1,\ldots,c'_t,c'_{t+1}}(G^*)\cap\gamma_{c'_{t+1}+1}(\gamma_{c'_t+1}(\cdots(\gamma_{c'_1+1}(G^*)))\cdots))$.\\
We consider the following notations, inductively:\\
$\rho_0(G^*)=G^*$ and $\rho_i(G^*)=\gamma_{c_i+1}(\rho_{i-1}(G^*)),$ for all $i\geq 1$,\\
$\rho'_0(G^*)=G^*$ and $\rho'_i(G^*)=\gamma_{c'_i+1}(\rho'_{i-1}(G^*)),$ for all $i\geq 1$.\\
Since, $\rho_{t+1}(G)=1$, so we have $\rho_{t+1}(G^*/A)=1$, and hence
$\rho_{t+1}(G^*)\subseteq A$. Also $A\subseteq\rho'_{t+1}(G^*),$
then $\rho_{t+1}(G^*)\subseteq\rho'_{t+1}(G^*)$. On the other
hand, by $c'_i\geq c_i$ for all $1\leq i\leq t$ and
$c'_{t+1}>c_{t+1}$ we can imply that
$\rho'_j(G^*)\subseteq\rho_j(G^*)$, for all $1\leq j\leq t+1$. Therefore
$$\rho'_{t+1}(G^*)=\rho_{t+1}(G^*)\ \ (I).$$
Consider the following trivial inclusions:\\
$\rho'_{t+1}(G^*)=\gamma_{c'_{t+1}+1}(\rho'_t(G^*))\subseteq\gamma_{c'_{t+1}}(\rho'_t(G^*))\subseteq\gamma_{c'_{t+1}-1}(\rho'_t(G^*))\subseteq$
\begin{center}
$\cdots\subseteq\gamma_{c_{t+1}+1}(\rho'_t(G^*))\subseteq\gamma_{c_{t+1}+1}(\rho_t(G^*))=\rho_{t+1}(G^*).$
\end{center}
Thus by the equality $(I)$, we can conclude that
$$\gamma_{c'_{t+1}+1}(\rho'_t(G^*))=\gamma_{c_{t+1}+1}(\rho'_t(G^*))\ (II).$$
Since $\rho_{t+1}(G^*)\subseteq A\subseteq N^*_{c'_1,\ldots,c'_t,c'_{t+1}}(G^*)$, we have
$$[\cdots[[\rho_{t+1}(G^*),_{c'_1}\rho'_0(G^*)],_{c'_2}\rho'_1(G^*)],\cdots,_{c'_{t+1}}\rho'_t(G^*)]=1.$$
Clearly $\rho'_{t}(G^*)\subseteq\rho'_{i}(G^*)$ for all
$0\leq i\leq t$, so by $(II)$, we can conclude that
$$[\cdots[[\gamma_{c_{t+1}+1}(\rho'_t(G^*)),_{c'_1}\rho'_t(G^*)],_{c'_2}\rho'_t(G^*)],\cdots,_{c'_{t+1}}\rho'_t(G^*)]=1.$$
and then
$\gamma_{c_{t+1}+1+c'_1+\cdots+c'_{t+1}}(\rho'_t(G^*))=1$. Put
$c=c_{t+1}+1+c'_1+\cdots+c'_{t+1}$ and $k=c'_{t+1}-c_{t+1}$. By
division algorithm, there are $q,r\in\bf Z$ such that $c=kq+r$,
where $r<k$. Put $j=\min\{i\in{\bf N} |ki+r\geq c'_{t+1}+1\}$. Then $kj+r\geq c'_{t+1}+1$ and $k(j-1)+r<c'_{t+1}+1$. Now, using $(II)$ we have
$1=\gamma_c(\rho'_t(G^*))=[\gamma_{c'_{t+1}+1}(\rho'_t(G^*)),_{c-c'_{t+1}-1}\rho'_t(G^*)]$\\
$=[\gamma_{c_{t+1}+1}(\rho'_t(G^*)),_{c-c'_{t+1}-1}\rho'_t(G^*)]$\\
$=\gamma_{c-k}(\rho'_t(G^*))$\\
$\vdots$\\
$=\gamma_{c-k(q-j)}(\rho'_t(G^*))$\\
$=\gamma_{kj+r}(\rho'_t(G^*))$\\
$=\gamma_{k(j-1)+r}(\rho'_t(G^*))$\\
$\supseteq\gamma_{c'_{t+1}+1}(\rho'_t(G^*))$\\
$=\rho'_{t+1}(G^*)$.
Hence $\rho'_{t+1}(G^*)=1$ and
so $A=1$, which is a contradiction. $\Box$\\
{\bf Notes}

(i)\ The condition $c'_{t+1}>c_{t+1}$ in the theorem 2.3 is
essential, since the first author [10] showed that for any natural
number $n$, there exists a nilpotent group $G$ of class $n$ such
that ${\cal N}_cM(G)\neq 1$ and $G$ has at least one ${\cal
N}_c$-covering group for all $c\leq n$.

(ii)\ In a joint paper with the first author [11], it is shown
that a finitely generated abelian group $G\cong{\bf
Z}_{n_1}\oplus{\bf Z}_{n_2}\oplus\ldots\oplus{\bf Z}_{n_k}$, where
$n_{i+1}|n_i$ for all $1\leq i\leq k-1$, has a nontrivial
polynilpotent multiplier, ${\cal N}_{c_1,\cdots,c_t}M(G)$, if
$k\geq 3$. Hence we can find many groups satisfying in conditions
of Theorems 2.1 and 2.3.

\end{document}